\documentclass[letterpaper, 10pt, conference]{ieeeconf}

\IEEEoverridecommandlockouts                              
\usepackage[cmex10]{amsmath}
\usepackage{cite,url}
\usepackage{bm}
\usepackage{amsmath,amssymb}
\usepackage{booktabs}
\usepackage{csvsimple}
\usepackage{color}
\usepackage{graphicx}

\newtheorem{theorem}{Theorem}

\title{\LARGE \bf
Robust Gas Pipeline Network Expansion Planning to \\Support Power System Reliability
}

\author{Kaarthik Sundar$^\star$, Sidhant Misra$^\dagger$, Anatoly Zlotnik$^\dagger$, Russell Bent$^\dagger$
\thanks{$^\star$\texttt{kaarthik@lanl.gov}, Information Systems \& Modeling, Los Alamos National Laboratory, Los Alamos, NM 87545}
\thanks{$^\dagger$\texttt{\{azlotnik,sidhant,rbent\}@lanl.gov}, Applied Mathematics \& Plasma Physics, Los Alamos National Laboratory, Los Alamos, NM 87545}}

\begin{document}
\maketitle

\begin{abstract}
We examine the problem of optimal transport capacity expansion planning for a gas pipeline network to service the growing demand of gas-fired power plants that are increasingly used to provide base load, flexibility, and reserve generation for bulk electric system. The aim is to determine the minimal cost set of additional pipes and gas compressors that can be added to the network to provide the additional capacity to service future loads. This combinatorial optimization problem is initially formulated as a mixed-integer nonlinear program, which we then extend to account for the variability that is inherent to the demands of gas-fired electricity production and uncertainty in expected future loads.  We consider here steady-state flow modeling while ensuring that the solution is feasible for all possible values of interval uncertainty in loads, which results in a challenging semi-infinite problem.  We apply previously derived monotonicity properties that enable simplification of the problem to require constraint satisfaction in the two extremal scenarios only, and then formulate the robust gas pipeline network expansion planning problem using a mixed-integer second order cone formulation.  We consider case studies on the Belgian network test case to examine the performance of the proposed approach. 
\end{abstract}


\section{Introduction}
The imperative to lower emissions associated with electric power production, as well as an increased supply and lower prevailing prices of natural gas, have led to widespread construction of gas-fired electric generators throughout the world \cite{Lyons2013}. Such generators now supply bulk electric systems with base load, spinning reserves, and balancing services that compensate for the variability of intermittent renewable energy sources such as wind and solar \cite{Rinaldi2001,Li2008}. Both use cases result in significant and (sometimes) rapid variation in natural gas consumption, which is increasingly difficult for operators of pipeline systems to service responsively. Extreme events, such as polar vortex conditions (extreme cold temperatures) that occur in the Northeast United States, further push interdependent electricity production and gas transportation infrastructure systems to their limits \cite{babula2014cold}. In these events, icy and windy weather strain generating stations, power lines, and other interdependent infrastructure systems that deliver that power, at the same time as cold temperatures induce increased demand for commercial and residential heating, so that the pipeline transportation capacity reserved by local utilities my become fully subscribed.  Because many power systems now rely on natural gas to supply base demand throughout the year, it is crucially important to design and expand natural gas pipeline networks to handle extreme events and meet increasing demands driven by the need to provide fuel for electric power generation.  However,   capacity expansion for a natural gas pipeline requires substantial capital expenditures, so it is critical for these resources to be allocated most effectively to cover a wide range of load scenarios that may arise from variations in weather, climate, power system function, and future resource allocation \cite{qiu2014low,qiu2014multi}. 

In a recent study \cite{Borraz2016}, the authors consider a network expansion planning problem for natural gas pipelines where the cost of infrastructure that is able to support fixed loads is minimized.  While these results provide a conceptual foundation as well as the modeling and computational basis for optimizing gas pipeline capacity expansion, a problem that considers only a single future scenario may yield conservative or insufficient results \cite{sanchez2016convex}.  Uncertain factors that may lead to sub-optimal decisions are related to resource allocation, economics, and climate change \cite{chaudry2014combined}.  Furthermore, even if parameters related to these factors are predicted exactly, the loads on a gas pipeline network would still have large variations caused by the need to service the bulk electric system.

In this study we address the issue of uncertainty in pipeline load scenarios that arises from seasonal and daily variations in bulk electric system load, assuming no uncertainty in the basic structure of the system and baseline loads.  Thus, we develop a systematic method for constructing a meaningful representation to capture the variety of possible load scenarios for a gas pipeline network by considering gas-fired generator fuel burn profiles and heat-rate curves. Specifically, we construct gas withdrawal uncertainty sets and use these to formulate a robust gas expansion planning problem (RGEPP) formulation for natural gas delivery systems.  The RGEPP is formulated first as a semi-infinite mixed-integer nonlinear program (MINLP), where loads at pipeline network nodes with gas-fired generators may take one of a continuum of values.  

A major advance in this current study beyond previously published work is the application of monotone system theory \cite{hirsch2006monotone} in a manner that has specific advantages for controlled natural gas flow networks, which enables a significantly simplified expansion planning problem representation and promotes desirable solution properties for combinatorial optimization solvers.  Monotone system theory has previously been applied to the case of natural gas flow networks by the authors and collaborators \cite{dvijotham2015natural,Zlotnik2016}, and was investigated for nonlinear programming applications \cite{vuffray2015monotonicity}.  Here we apply a similar insight to re-formulate the semi-infinite RGEPP problem into a MINLP with a finite number of specific constraints, in which the range of uncertainty is accounted for using a specifically defined scenarios.  This in turn is reformulated using McCormick relaxations \cite{mccormick1976} into a convex mixed integer second order cone problem (MISOCP), for which solutions can be readily computed by standard commercially available solvers. To evaluate the advantages of the proposed formulation, we perform extensive computational experiments on the widely-examined ``Belgium'' gas pipeline test network.  We benchmark computation time and perform comparisons with random sampling of the uncertainty space in order to demonstrate the efficiency of our approach, and also examine how the further removal of certain constraints based on the monotone order principle affect the resulting expansion solutions.  The key property of the robust solution is that \emph{the resulting pipeline system expansion plan will support any bulk electric system dispatch schedule that lies within a defined uncertainty set.}

The rest of the manuscript is structured as follows.  We formally state the problem, notation, and full formulation in Section \ref{sec:statement}, and then describe the reformulation procedure in Section \ref{sec:reformulation}.  The application of monotonicity theorems to the reformulated problem is described in Section \ref{sec:monotonicity}.  The computational studies and implementation details are presented in Section \ref{sec:results}, and we conclude with Section \ref{sec:conclusion}.



\section{Problem Statement and Formulation} \label{sec:statement}
The goal of the robust gas network expansion planning problem (RGEPP) is to determine the minimum cost set of capital expansions for a gas pipeline network to increase its capacity to support projected increases in demand, which are subject to interval uncertainty. Formally, the RGEPP is stated as follows: given a natural gas pipeline network with nodes $\mathcal N$, pipelines, $\mathcal P_e$, and gas compressor stations $\mathcal C_e$, and potential expansion options for pipelines from the set $\mathcal P_n$ and compressors in the set $\mathcal C_n$, compute a minimum cost network expansion plan for to meet the increased and uncertain gas demand while satisfying the flow and pressure limits of the pipeline system. The problem is initially formulated as a Mixed-Integer Non-Linear Program (MINLP), and a Mixed-Integer Second Order Cone (MISOC) relaxation is developed and presented in the next section. First, we present the notation that will be used throughout the rest of the article.


\subsection{Notation} \label{sec:notations}

\noindent \textit{Sets:} \\
 \noindent $\mathcal N$, $\mathcal C_e$, $\mathcal P_e$ - (existing) nodes, compressors, pipes \\ 
 \noindent $\mathcal C_n$, $\mathcal P_n$ - (potential new) compressors, pipes \\
 \noindent $\mathcal C = \mathcal C_e \cup \mathcal C_n$, $\mathcal P = \mathcal P_e \cup \mathcal P_n$ \\
 \noindent $\mathcal{G}, \mathcal R$, $\mathcal D$ - generation, receipt and delivery points \\ 
 \noindent $\mathcal R(i)$, $\mathcal D(i)$ - receipt and delivery points at node $i$ \\  
 \noindent $\mathcal K$ - scenario profiles \\
 \noindent $\mathcal{S}_k$ - set of scenarios in scenario profile $k \in \mathcal{K}$ \\
 \noindent $\mathcal S = \cup_{k \in \mathcal{K}} \mathcal{S}_k$ - set of all demand scenarios \\
 \noindent $\mathcal E(i)$ - subset of pipes and compressors (existing and new) connected to node $i$ and oriented \textit{from} $i$ \\ 
 \noindent $\mathcal E^r(i)$ - subset of pipes and compressors (existing and new) connected to node $i$ and oriented \textit{to} $i$ \\ 
 \noindent Note: Compressors and pipes are denoted by triplets $(e,i,j)$ consisting of a unique identifier $e$ linking nodes $i$ and $j$. For convenience, such a triplet $(e,i,j)$ will be denoted by $e_{ij}$ or $e$, interchangeably. \\
 \noindent \textit{Decision variables:} \\ 
 \noindent $p_i^s$ - pressure at node $i$ for $s \in \mathcal S$ \\
 \noindent $\pi_i^s$ - square of pressure at node $i$ for $s \in \mathcal S$ \\ 
 \noindent $f_e^s$ - mass flow rate across $e \in \mathcal C \cup \mathcal P$  for $s \in \mathcal S$\\
 \noindent $s_i^s$ - total gas produced at receipt points in $\mathcal R(i)$ for $s \in \mathcal S$ \\
 \noindent $\gamma_e^s$ - auxiliary variable for each pipe $e \in \mathcal P$ for $s \in \mathcal S$ \\
 \noindent $y_e^s$ - binary flow direction variable for $e \in \mathcal C \cup \mathcal P$ for $s \in \mathcal S$\\
 \noindent $z_e$ - binary expansion variable for each $e \in \mathcal C_n \cup \mathcal P_n$ \\
 \noindent \textit{Parameters:} \\ 
 \noindent $d_i^s$ - gas demand at a node $i\in\mathcal N$ in scenario $s\in\mathcal S$ \\
 \noindent $[{\bm d}^{\ell,k}_i$, ${\bm d}^{u,k}_i]$ - uncertainty set of total gas demand for all delivery points in $\mathcal D(i)$ for a profile $k\in\mathcal K$ \\
 \noindent $\bm w_e$ - resistance of the pipe $e \in \mathcal P$ \\ 
 \noindent $\bm a$ - speed of sound \\ 
 \noindent $\bm \beta_e$ - friction factor of the pipe $e \in \mathcal{P}$ \\ 
 \noindent $\bm{\ell}_e$, $\bm D_e$ - length, diameter of the pipe $e \in \mathcal{P}$  \\ 
 \noindent [$\underline{\bm \pi}_i$, $\overline{\bm \pi}_i$] - min and max limits for $\pi_i$ \\ 
 \noindent [$\underline{\bm \alpha}_e$, $\overline{\bm \alpha}_e$] - min and max compression ratio limits for $e \in \mathcal C$ \\
 \noindent $\bm f_e$ - max flow rate  for $e \in \mathcal C \cup \mathcal P$ \\
 \noindent $\bm c_e$ - cost of building component $(e, i, j) \in \mathcal C_n \cup \mathcal P_n$ \\
 \noindent Note: A scenario $\bm s\in\mathcal S_k$ is characterized by the specific load values $d_i^s\in [{\bm d}^{\ell,k}_i$, ${\bm d}^{u,k}_i]$ for  $i\in\mathcal N$. A profile of demand scenarios $\mathcal S_k$ is the Cartesian product of all the uncertain gas demand values, i.e., $\mathcal S_k=\prod_{i\in\mathcal N}[{\bm d}^{\ell,k}_i$, ${\bm d}^{u,k}_i]$. With multiple scenario profiles, the total uncertainty set of possible scenarios is given as $\mathcal{S} = \bigcup_{k\in \mathcal K} \mathcal S_k$, where $\mathcal K$ is the collection of profiles. In the case of only one profile $\mathcal S$, we may ignore the index $k$.

\subsection{Steady State Gas Flow Physics}
Here we review the physical laws that govern steady flow of natural gas through pipelines. 
The physics of flow across a  pipeline, $e = (i,j)$, are  described  by  a  set  of  partial  differential equations (PDEs) in time and space \cite{Thorley1987}. 
In the steady-state, the PDEs reduce to equations of the form
\begin{flalign}
\quad \pi_i - \pi_j = \bm w_e f_e | f_e |, \label{eq:steady-state-physics}
\end{flalign}
where the phenomenological expression on the right hand side quantifies the dissipation of kinetic energy that results from eddies formed by turbulent flow through the pipe.  The parameter $\bm w_e$ is called a resistance factor, and is given by 
\begin{flalign}
\bm w_e = \frac{4 \bm \beta_e \bm l_e \bm a^2}{\pi^2 \bm D_e^5 }. \label{eq:resistance}
\end{flalign}
For a detailed derivation of this equation, interested readers are referred to \cite{Sundar2018}.  To compensate for the dissipation of energy along the direction of flow, a gas pipeline utilizes compressors to boost flow and pressure throughout the system.  We model these components as short pipes with zero resistance values, which create a jump in pressure while preserving flow in the direction of the compressor's orientation. When the gas flows through the compressor in the opposite direction of its orientation, the compressor is assumed to not offer any pressure boost. The elements $(e, i, j) \in \mathcal C$ may be used to represent compressor stations or individual compressor machines, and are modeled as short pipes with zero resistance values that can increase pressure according to $\pi_j^s=\alpha_e^2\pi_i^s$ in the direction of the compressor's orientation i.e., from $i$ to $j$. 

\subsection{Problem Formulation} \label{sec:formulation}
Given the notations in Sec. \ref{sec:notations}, the RGEPP is formulated as follows:
\begin{subequations}
\begin{flalign}
& \min \quad \sum_{(p,i,j) \in \mathcal P_n} \bm c_{p} z_p + \sum_{(c,i,j) \in \mathcal C_n} \bm c_c z_c & \label{eq:obj} \\ 
& \mbox{subject to} \\
& \text{for every scenario $s \in \mathcal S$:} \notag \\ 
& \pi_i^s - \pi_j^s = \bm w_e |f_e^s| f_e^s \quad \forall (e, i, j) \in \mathcal P_e, & \label{eq:pipe-physics-existing} \\ 
& \pi_i^s - \pi_j^s = 0, \text{ if $f_e^s \leqslant 0$, }  \forall (e, i, j) \in \mathcal C_e, & \label{eq:compressor-physics-1-existing} \\ 
& \pi_j^s \in [\underline{\bm \alpha}_e^2 \pi_i^s, \overline{\bm \alpha}_e^2 \pi_i^s], \text{if $f_e^s \geqslant 0$,}  ~\forall (e, i, j) \in \mathcal C_e, & \label{eq:compressor-physics-2-existing} \\ 
& \sum_{e \in \mathcal E(i)} f_e^s - \sum_{e \in \mathcal E^r(i)} f_e^s = s_i^s - d_i^s ~ \forall i \in \mathcal N,  & \label{eq:flow-balance} \\ 
& \pi_i^s - \pi_j^s = \bm w_e |f_e^s| f_e^s \text{ if $z_e = 1$, } \forall (e, i, j) \in \mathcal P_n, & \label{eq:pipe-physics-new} \\ 
& \pi_i^s - \pi_j^s = 0, \text{ if $f_e^s \leqslant 0$, $z_e = 1$, }  \forall (e, i, j) \in \mathcal C_n, & \label{eq:compressor-physics-1-new} \\ 
& \pi_j^s \in [\underline{\bm \alpha}_e^2 \pi_i^s, \overline{\bm \alpha}_e^2 \pi_i^s], \text{if $f_e^s \geqslant 0$, $z_e = 1$} ~\forall (e, i, j) \in \mathcal C_n, & \label{eq:compressor-physics-2-new} \\
& \underline{\bm \pi}_i \leqslant \pi_i^s \leqslant \overline{\bm \pi}_i \quad \forall i \in \mathcal N, & \label{eq:pressure-bounds} \\ 
& -{\bm f}_e \leqslant f_e^s \leqslant {\bm f}_e \quad \forall e \in \mathcal C \cup \mathcal P. & \label{eq:flow-bounds}
\end{flalign}
\label{eq:formulation}
\end{subequations}
The formulation for the RGEPP, as stated in Eq. \eqref{eq:formulation}, is a MINLP formulation with disjunctive constraints. The objective in Eq. \eqref{eq:obj} aims to minimize the network expansion cost. The constraints in Eq. \eqref{eq:pipe-physics-existing} enforces the steady-state gas flow physics for each existing pipe $(e, i, j) \in \mathcal P_e$ in the network. The Eqs. \eqref{eq:compressor-physics-1-existing} and \eqref{eq:compressor-physics-2-existing} enforce the boosting limits of the existing compressors $(e, i, j) \in \mathcal C_e$ when the flow in the compressor is directed along the orientation of the compressor. The Eq. \eqref{eq:flow-balance} enforces the nodal balance condition for each node in the system. As for the constraints in Eq. \eqref{eq:pipe-physics-new} -- \eqref{eq:compressor-physics-2-new}, they enforce the steady state physics and the boosting limits for the new pipes and new pipe compressors, respectively. They are enforced if and only if the corresponding pipes and compressors are constructed. Finally, the constraints in Eq. \eqref{eq:pressure-bounds} and \eqref{eq:flow-bounds} enforce the bounds on the nodal pressure and flows on the pipes and compressors. The constraint in Eq. \eqref{eq:pipe-physics-existing} is a nonlinear constraint while the Eq. \eqref{eq:pipe-physics-new} is a nonlinear disjunctive constraint that is enforced only when the corresponding value of $z_e = 1$. As with compressors, the constraints for both the existing and new compressors in Eq. \eqref{eq:compressor-physics-1-existing} -- \eqref{eq:compressor-physics-2-existing} and Eq. \eqref{eq:compressor-physics-1-new} -- \eqref{eq:compressor-physics-2-new}, respectively, are linear disjunctive constraints. In the next section, we reformulate the disjunctive constraints using binary variables, resulting in a MINLP reformulation of Eq. \eqref{eq:formulation}. Subsequently, we then develop a MISOC relaxation using the reformulated MINLP.

\section{Reformulation and Convex Relaxation} \label{sec:reformulation}
A MINLP reformulation of the nonlinear and linear disjunctive constraints in Eq. \eqref{eq:formulation} is done  using binary flow direction variables $y_e^s$ for each $e \in \mathcal C \cup \mathcal P$. Given a pipe or a compressor $(e,i,j) \in \mathcal C \cup \mathcal P$, $y_e^s$ takes a value $1$ if the mass flow $f_e^s \geqslant 0$ and $0$, otherwise. We remark that if the flow $f_e^s \leqslant 0$, then gas is flowing from the node $j$ to node $i$. 

\subsection{Reformulation of Eq. \eqref{eq:pipe-physics-existing}}
Given these notations, Eq. \eqref{eq:pipe-physics-existing} for any $(e, i,j) \in \mathcal P_e$, can equivalently be reformulated to
\begin{subequations}
\begin{flalign}
& \gamma_e^s = \bm w_e (f_e^s)^2 & \label{eq:pipe-existing-1} \\ 
& \gamma_e^s \geqslant \pi_j^s - \pi_i^s + 2 y_e^s (\underline{\bm \pi}_i - \overline{\bm \pi}_j) & \label{eq:mc-1} \\
& \gamma_e^s \geqslant \pi_i^s - \pi_j^s + 2 (y_e^s - 1) (\overline{\bm \pi}_i - \underline{\bm \pi}_j) & \label{eq:mc-2} \\
& \gamma_e^s \leqslant \pi_j^s - \pi_i^s + 2 y_e^s (\overline{\bm \pi}_i - \underline{\bm \pi}_j) & \label{eq:mc-3} \\
& \gamma_e^s \leqslant \pi_i^s - \pi_j^s + 2 (y_e^s - 1) (\underline{\bm \pi}_i - \overline{\bm \pi}_j) & \label{eq:mc-4} \\ 
& - {\bm f}_e (1-y_e^s) \leqslant f_e^s \leqslant {\bm f}_e y_e^s & \label{eq:flow-bounds-pipes-existing} 
\end{flalign}
\label{eq:reformulation-pipes-existing}
\end{subequations}
where, $\gamma_e^s$ is an auxiliary variable for pipe $(e, i, j) \in \mathcal P_e$. The Eq. \eqref{eq:mc-1} -- \eqref{eq:mc-4} are the McCormick envelopes \cite{mccormick1976} for the equation $\gamma_e^s = (2y_e^s - 1) (\pi_i^s - \pi_j^s)$. These envelopes result in an exact reformulation because it is the product of a variable that takes a value of one or negative one, $(2y_e^s - 1)$, with a continuous variable, $(\pi_i^s - \pi_j^s)$. The Eq. \eqref{eq:flow-bounds-pipes-existing} bounds the mass flow on the pipe using the flow direction variable $y_e^s$. The only nonlinear constraint in the reformulation is given by Eq. \eqref{eq:pipe-existing-1}. 

\subsection{Reformulation of Eq. \eqref{eq:pipe-physics-new}}
The constraint in Eq. \eqref{eq:pipe-physics-new} for every $(e, i, j) \in \mathcal P_n$ is similarly reformulated using binary flow direction variable $y_e^s$ as 
\begin{flalign}
& z_e \gamma_e^s = \bm w_e (f_e^s)^2 \text{ and Eq. \eqref{eq:mc-1} -- \eqref{eq:flow-bounds-pipes-existing}. } & \label{eq:reformulation-pipes-new} &
\end{flalign}
The Eq. \eqref{eq:reformulation-pipes-new} ensures that flow across a new pipeline $(e, i, j)$ is non-zero only if the pipeline is built by forcing the variable $z_e$ to take a value of $1$. 

\subsection{Reformulation of Eqs. \eqref{eq:compressor-physics-1-existing} and \eqref{eq:compressor-physics-2-existing}} 
The constraints enforced on the existing compressors in the network are reformulated into linear constraints using binary flow direction variables $y_e^s$ as shown below:
\begin{subequations}
\begin{flalign}
& y_e^s (\underline{\bm \pi}_i - \overline{\bm \pi}_j) \leqslant \pi_i^s - \pi_j^s \leqslant y_e^s (\overline{\bm \pi}_i - \underline{\bm \pi}_j) & \label{eq:compressor-existing-reformulation-1} \\ 
& \underline{\bm \alpha}_e^2 \pi_i^s + (1-y_e^s) (\underline{\bm \pi}_j -\underline{\bm \alpha}_e^2 \overline{\bm \pi}_i) \leqslant \pi_j^s & \label{eq:compressor-existing-reformulation-2} \\
& \pi_j^s \leqslant \overline{\bm \alpha}_e^2 \pi_i^s + (1-y_e^s) (\overline{\bm \pi}_j - \overline{\bm \alpha}_e^2 \underline{\bm \pi}_i) & \label{eq:compressor-existing-reformulation-3} \\
& - {\bm f}_e (1-y_e^s) \leqslant f_e^s \leqslant {\bm f}_e y_e^s & \label{eq:flow-bounds-compressors} 
\end{flalign}
\label{eq:compressor-existing-reformulation} 
\end{subequations}
The constraints in Eq. \eqref{eq:compressor-existing-reformulation-1} and \eqref{eq:compressor-existing-reformulation-2} -- \eqref{eq:compressor-existing-reformulation-3} ensure Eq. \eqref{eq:compressor-physics-1-existing} and respectively, and \eqref{eq:compressor-physics-2-existing} are satisfied when the gas flow in the compressor is in the reverse direction and forward direction, respectively.

\subsection{Reformulation of Eqs. \eqref{eq:compressor-physics-1-new} and \eqref{eq:compressor-physics-2-new}}
The constraints for each new compressor $(e, i, j) \in \mathcal C_n$ are reformulated into linear constraints using the binary flow direction variables $y_e^s$ and expansion variables $z_e$ as follows:
\begin{subequations}
\begin{flalign}
& (2 + y_e^s - z_e) (\underline{\bm \pi}_i - \overline{\bm \pi}_j) \leqslant \pi_i^s - \pi_j^s & \label{eq:compressor-new-reformulation-1} \\ 
& \pi_i^s - \pi_j^s \leqslant (2 + y_e^s - z_e) (\overline{\bm \pi}_i - \underline{\bm \pi}_j) & \label{eq:compressor-new-reformulation-3} \\
& \underline{\bm \alpha}_e^2 \pi_i^s + (2 - y_e^s - z_e) (\underline{\bm \pi}_j -\underline{\bm \alpha}_e^2 \overline{\bm \pi}_i) \leqslant \pi_j^s & \label{eq:compressor-new-reformulation-2} \\
& \pi_j^s \leqslant \overline{\bm \alpha}_e^2 \pi_i^s + (2 - y_e^s - z_e) (\overline{\bm \pi}_j - \overline{\bm \alpha}_e^2 \underline{\bm \pi}_i) & \label{eq:compressor-new-reformulation-4} \\
& - {\bm f}_e (1-y_e^s) \leqslant f_e^s \leqslant {\bm f}_e y_e^s & \label{eq:flow-bounds-compressors-new} \\
& - {\bm f}_e z_e \leqslant f_e^s \leqslant {\bm f}_e z_e & \label{eq:flow-bounds-compressors-expansion} 
\end{flalign}
\label{eq:compressor-new-reformulation}
\end{subequations}
Similar to the previous section, the Eq. \eqref{eq:compressor-new-reformulation-1} and \eqref{eq:compressor-new-reformulation-3} are equivalent to \eqref{eq:compressor-physics-1-new} and Eq. \eqref{eq:compressor-new-reformulation-2} and \eqref{eq:compressor-new-reformulation-4} are equivalent to \eqref{eq:compressor-physics-2-new}. As for the constraints in Eq. \eqref{eq:flow-bounds-compressors-new} and \eqref{eq:flow-bounds-compressors-expansion}, they restrict the gas flow and its direction based on the values taken by $y_e^s$ and $z_e$. 

In summary, the MINLP reformulation for the RGEPP is given by 
\begin{flalign*}
& \min \quad \sum_{(p,i,j) \in \mathcal P_n} \bm c_{p} z_p + \sum_{(c,i,j) \in \mathcal C_n} \bm c_c z_c & \\
& \text{subject to, for every scenario $s \in \mathcal S_k$, for all $k\in\mathcal K$:} \notag \\  
& \qquad \text{\eqref{eq:reformulation-pipes-existing}, \eqref{eq:reformulation-pipes-new}, \eqref{eq:compressor-existing-reformulation}, \eqref{eq:compressor-new-reformulation}, \eqref{eq:flow-balance}, \eqref{eq:pressure-bounds}, and \eqref{eq:flow-bounds}.}
\end{flalign*}

\subsection{MISOC relaxation} \label{subsec:misoc}
The MINLP reformulation for the RGEPP still a semi-infinite problem because the number of scenarios are infinite and to the best of our knowledge there is no algorithm in the literature that can directly solve semi-infinite MINLPs to global optimality. Furthermore, even if there are only a finite number of scenarios, developing an algorithm that can solve the resulting MINLP to global optimality is nonetheless challenging \cite{oertel2014complexity}. With this in mind, we also develop an MISOC relaxation of the MINLP that can be solved to global optimality using constraint generation techniques. The only nonlinear constraints in the MINLP are Eq. \eqref{eq:pipe-existing-1} and \eqref{eq:pipe-physics-new}, and both of these constraints can be relaxed to second-order conic (SOC) and rotated second-order conic (RSOC) constraints, which then take the forms
\begin{flalign}
& \gamma_e^s \geqslant \bm w_e (f_e^s)^2 \quad \forall (e, i, j) \in \mathcal P_e, & \label{eq:soc-existing} \\ 
& z_e \gamma_e^2 \geqslant \bm w_e (f_e^s)^2 \quad \forall (e, i, j) \in \mathcal P_n. & \label{eq:rsoc-new}
\end{flalign}
Equations \eqref{eq:soc-existing} and \eqref{eq:rsoc-new} are SOC and RSOC relaxations of equations \eqref{eq:pipe-existing-1} and the equality in \eqref{eq:reformulation-pipes-new}, respectively. Hence, the MISOC relaxation of the RGEPP is obtained by replacing Eq. \eqref{eq:pipe-existing-1} and the equality in \eqref{eq:reformulation-pipes-new} in the MINLP reformulation with equations \eqref{eq:soc-existing} and \eqref{eq:rsoc-new}, respectively.  This relaxation approach facilitates computational solution of the finite-scenario network expansion problem.  In the next section, we present monotonicity results for gas pipeline networks where flows are actuated by controllable compressors.  These results are applied below to reduce the semi-infinite MISOC relaxation to a finite MISOC where a collection of extremal scenarios are used to characterize the uncertainty set.

\section{Monotonicity in Gas Pipeline Flows} \label{sec:monotonicity}
In both the steady-state and transient regimes, the natural gas pipeline flow physics possess certain monotonicity properties \cite{vuffray2015monotonicity,Zlotnik2016,misra2016monotone}. Essentially, these properties state that, all else held equal, the pressures observed in a high demand scenario are lower than the pressures in a low demand scenario.  Application of these properties has been shown to greatly simplify the semi-infinite optimization problem corresponding to optimally controlling the operations of gas networks under uncertain demands. Indeed, as shown in e.g. \cite{vuffray2015monotonicity}, when the uncertainty set is described by a box, the semi-infinite problem can be reduced to a two-scenario problem where all demands are set to the maximum and minimum possible value respectively. 

However, the application of the monotonicity result in the preceding study on gas pipeline network optimization \cite{vuffray2015monotonicity} assumes that the compression ratios are kept constant over different realizations of the uncertainty, which is a reasonable assumption when considering short-term intra-day operations, but is an unrealistic assumption for the expansion planning problem that is considered here. The time scale of the expansion planning problem is much larger, because the scenario set must account for the possible operating conditions that could be encountered throughout several months or even multiple years. Once the expansion plan is implemented and a certain operating scenario is realized, the compressors are expected to be adjusted to the realized scenario in order to ensure feasibility of intra-day operations. However, in the following sections we propose that a generalized version of the monotonicity principle that allows for compressors to react to the uncertainty realization can be used to simplify the semi-infinite MISOC described in Section \ref{sec:reformulation} into a standard MISOC. 

\subsection{Monotonicity of Gas Flows Under Non-Decreasing Compressor Recourse}
A generalization of the monotonicity property \emph{with} recourse is established in \cite{misra2020monotone}. We state the theorem for completeness.
\begin{theorem} \label{thm:monotonicity}
    Let $\mathcal{N}_1 \subseteq \mathcal{N}$ be the set of nodes without withdrawals, and $\mathcal{N}_2  = \mathcal{N} \setminus \mathcal{N}_1$ be the set of nodes with active demands. Consider two different demand scenarios $\bm{d}^{(low)}$ and $\bm{d}^{(high)}$ such that for all $i \in \mathcal{N}_2$ we have $\bm{d}_i^{(low)} \leqslant \bm{d}_i^{(high)}$, and for all $i \in \mathcal{N}_1$, we have $\pi_i^{(low)} \geqslant \pi_i^{(high)}$.
    Additionally assume that for each $(e,i,j)$ with compressors, the output pressure $\pi_j$ is a monotonically increasing function of the input pressure $\pi_i$.
    Then the nodal pressures in the system satisfy the ordering $\pi_i^{(low)} \geqslant \pi_i^{(high)}$ for all $i \in \mathcal{N}_2$.
\end{theorem}

\subsection{Monotonic Compression Policy Constraints}
We propose to use Theorem~\ref{thm:monotonicity} to simplify the semi-infinite MISOCP to a problem that is constrained to satisfy only a finite number of scenarios. To satisfy the premises of Theorem~\ref{thm:monotonicity}, we need to define  a non-decreasing compressor policy for all possible demands \emph{within} each scenario profile $s \in \mathcal{S}_k$. 

For the compressor policy, we propose the following:
\begin{align}   \label{eq:compression_policy}
    \pi_j^s - \pi_i^s = \eta_e^s, \quad \forall (e,i,j) \in \mathcal{C},
\end{align}
where $\eta_e^s \geqslant 0$ is a scenario $s$-dependent non-negative \emph{variable}. 
The proposed policy has the following desirable properties:
\subsubsection{Non-decreasing nature and scenario reduction}
The compression policy in \eqref{eq:compression_policy} is clearly non-decreasing and the outlet pressure increases with the input pressure. Therefore, for a given scenario $s\in \mathcal{S}_k$, the assumptions of Theorem~\ref{thm:monotonicity} are satisfied, and the semi-infinite program in Section \ref{subsec:misoc} can be reduced to a two-scenario program, for which the full reformulation is given in Section~\ref{subsec:two_scenario_problem}.
\subsubsection{Simple algebraic structure} The policy can be enforced in the expansion planning formulation with the linear constraint in \eqref{eq:compression_policy}, thus not adding any significant complexity to the MISOC formulation.
\subsubsection{Inner approximation} Fixing the compressor control policy according to \eqref{eq:compression_policy} can be viewed as an \textit{inner approximation} to the constraint sets \eqref{eq:compressor-existing-reformulation} and \eqref{eq:compressor-new-reformulation}.  Therefore, if the problem admits a feasible solution with the proposed compression policy, then it is guaranteed to admit a feasible solution with a general compression policy.
\subsubsection{Intuitive direction of adjustment}
Because the resulting formulation is an inner approximation, the quality of the solution depends on how restrictive the proposed policy is. An overly restrictive policy can result in highly suboptimal solutions with much higher costs. First, the quantity $\eta_e^s$ is an optimization variable that is dependent on the scenario $s \in \mathcal{S}_k$, thus allowing it to be chosen optimally during the solution process. Second, the policy adjusts the compression ratio $\bm \alpha_e$ in a natural way - it uses a higher compression ratio when the demand are higher and we expect the input pressure to the compressor to be lower. Indeed, assuming that the pressure at a defined-pressure (slack) node in the network is fixed, then for two possible demand realizations of the demand corresponding to low and high load, we expect that the input pressure in the low load case denoted by $\pi_i^{(high)}$ is higher than the input pressure corresponding to the high demand case $\pi_i^{(low)}$, i.e., $\pi_i^{(low)} > \pi_i^{(high)}$. Then, a simple computation shows that the compression ratio in the low load case is lower than the compression ratio in the high load case:
\begin{align}   \label{eq:compression_ratio_comparision}
    \bm \alpha_e^{(low)} = \left(\frac{\pi_i^{(low)} + \eta_e^s }{ \pi_i^{(low)}}\right)^{0.5} < \notag \\ \left(\frac{\pi_i^{(high)} + \eta_e^s }{ \pi_i^{(high)}}\right)^{0.5} = \bm \alpha_e^{(high)}, 
\end{align}
where the inequality follows by using $\pi_i^{(low)} > \pi_i^{(high)}$. This natural property of the policy in \eqref{eq:compression_policy} makes it mimic the optimal oracle policy better, thus reducing optimality gap.

\subsection{Two-Scenario Problem}   \label{subsec:two_scenario_problem}
We can now augment the MISOC problem developed in Section \ref{sec:reformulation} with the monotonic policy developed above.  We first define the MISOC constraint set for a scenario $\bm s\in\mathcal S$ by $\bm \Gamma(\bm s)$, which is given as
$$
\bm \Gamma(\bm s) \,\, \left\{\begin{array}{l}
\eqref{eq:soc-existing} \text{ and } \eqref{eq:mc-1}-\eqref{eq:flow-bounds-pipes-existing}, \quad \forall (e, i, j) \in \mathcal P_e, \\
\eqref{eq:rsoc-new} \text{ and } \eqref{eq:mc-1}-\eqref{eq:flow-bounds-pipes-existing}, \quad \forall (e, i, j) \in \mathcal P_n, \\
\eqref{eq:compression_policy} \text{ and } \eqref{eq:compressor-existing-reformulation}, \quad \forall (e,i,j) \in \mathcal C_e, \\
\eqref{eq:compression_policy} \text{ and } \eqref{eq:compressor-new-reformulation}, \quad \forall (e,i,j) \in \mathcal C_n, \\
\eqref{eq:flow-balance}, \eqref{eq:pressure-bounds}, \text{ and } \eqref{eq:flow-bounds}. \\
\end{array}\right.
$$
Let us denote by $\bm s_k^\ell$ and $s_k^u$ the scenarios corresponding to $d_i \equiv \bm d_i^{\ell,k}$ for all $i\in\mathcal N$ and $d_i \equiv \bm d_i^{u,k}$ for all $i\in\mathcal N$, respectively, for a profile $k\in\mathcal K$.   The final MISOCP reformulation of the RGEPP, including the monotonic compression policy constraints, and constant pressure constraints at gas production nodes, can be posed as 
\begin{flalign*}
 \min \quad &\sum_{(p,i,j) \in \mathcal P_n} \bm c_{p} z_p + \sum_{(c,i,j) \in \mathcal C_n} \bm c_c z_c  \\
 \text{subject to: }  & \forall k \in \mathcal{K} \begin{cases} \bm \Gamma(\bm s_k^\ell) \text{ and } \bm \Gamma(\bm s_k^u), \\
 \pi_i^{s_k^\ell} = \pi_i^{s_k^u} \ i \in \mathcal{G}.
 \end{cases}
\end{flalign*}
The above formulation is a computationally tractable representation of the RGEPP.  We examine the computational properties of our implementation in the following section.

\section{Computational Results} \label{sec:results}
We present the results of extensive computational studies to show the effectiveness of the mathematical programming formulations proposed in Section \ref{sec:reformulation} in conjunction with application of the monotonicity results in Section \ref{sec:monotonicity}. We first present the benchmark test case, and an overview of the computational experiments that are performed on the benchmark instances.  
\subsection{Description of Benchmarks} \label{subsec:case-studies}
\begin{figure}
    \hspace{-1.2cm}\includegraphics[width=1.25\linewidth]{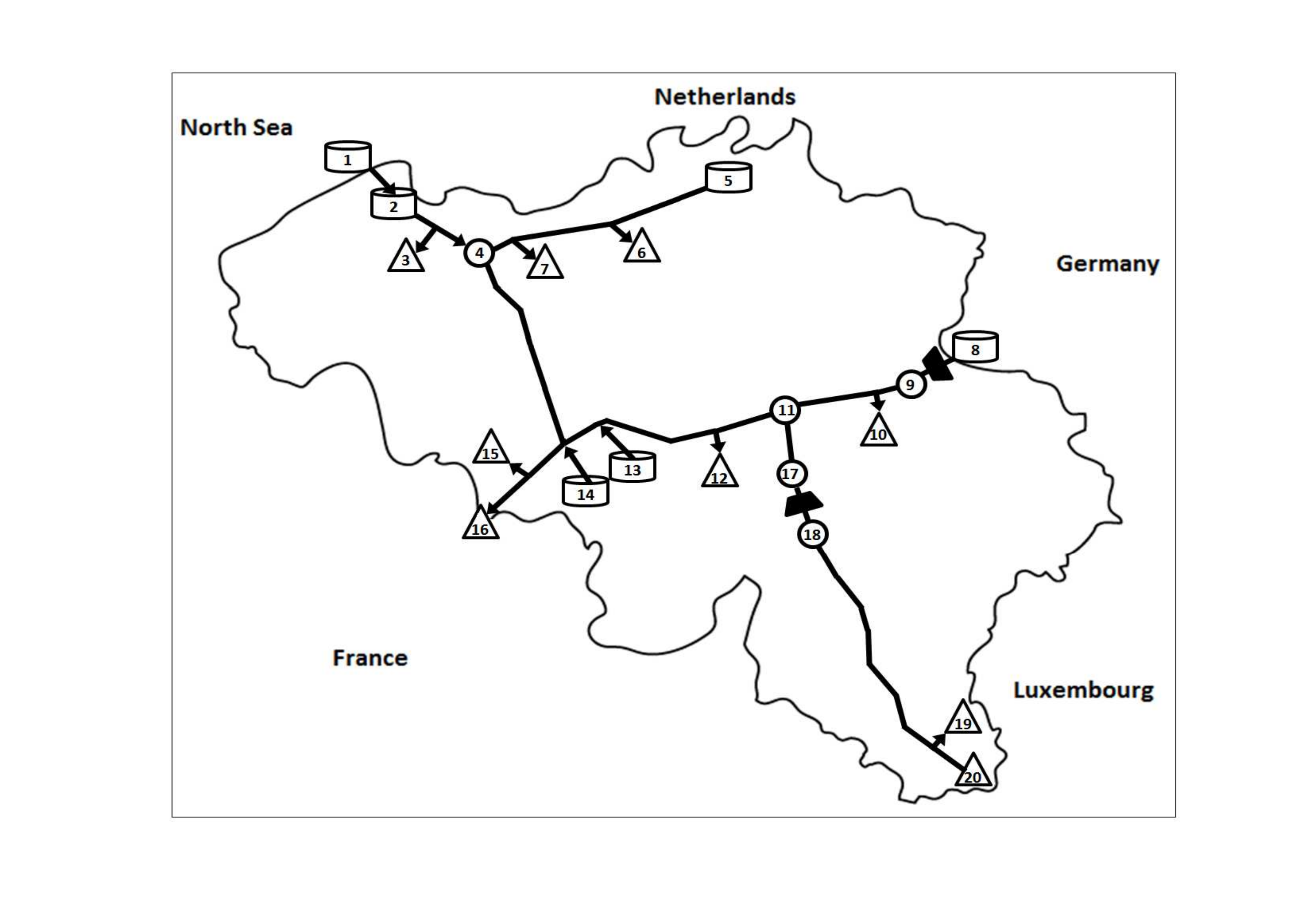}
    \vspace{-9ex}
    \caption{Belgian gas network schematic}
    \label{fig:belgian}
\end{figure}

All the computational experiments were performed on the Belgian network \cite{De1996optimal} shown in Fig. \ref{fig:belgian}. Three different benchmarks were created from the base Belgian network with different number of nodes and different number of new expansion options that can potentially be added to the existing network \cite{Borraz2016}. The overall size parameters of the nominal network model and the three expansion test cases are shown in Table \ref{tab:belgian-info}.

\begin{figure*}
    \centering
    \includegraphics[width=.95\linewidth]{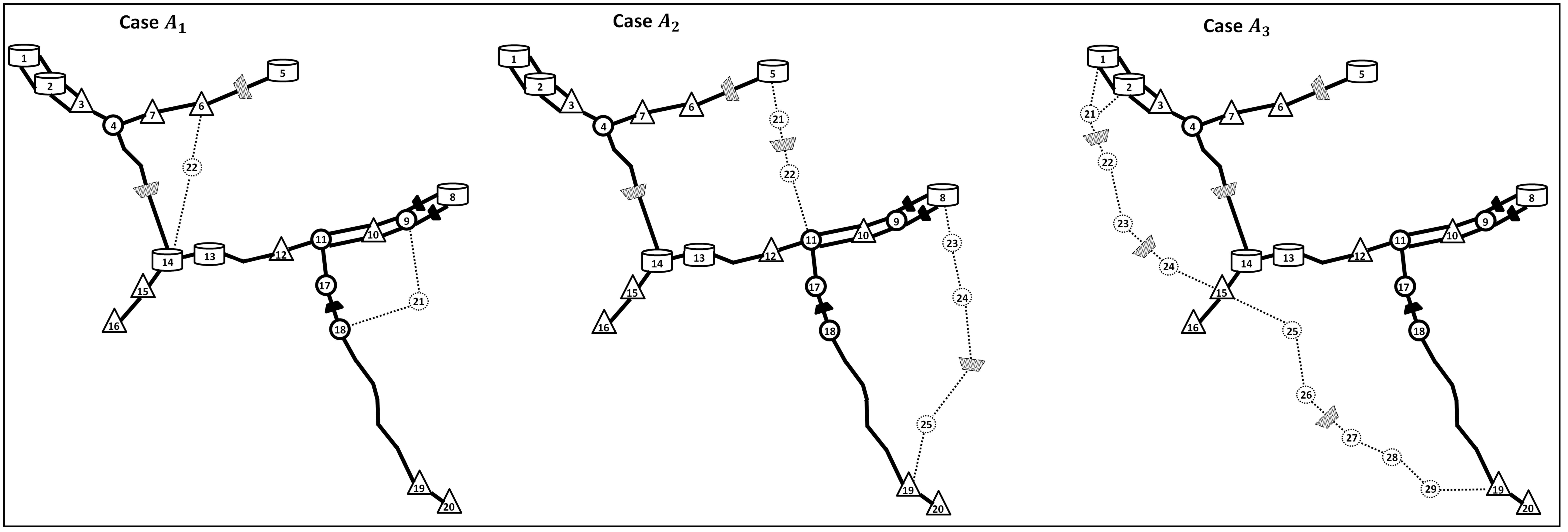}
    \caption{Schematic for expansion benchmarks $A_1$, $A_2$, and $A_3$.}
    \label{fig:expansion-belgian}
        \vspace{-2ex}
\end{figure*}

\begin{table}
    \centering
    \vspace{-2ex}
    \caption{The Belgian network benchmarks}
    \label{tab:belgian-info}
    \begin{tabular}{c|ccccc}
    \toprule 
        Network name & $|\mathcal N|$ & $|\mathcal P_e|$ & $|\mathcal C_e|$ & $|\mathcal P_n|$ & $|\mathcal C_n|$ \\
    \midrule 
        Base & 20 & 24 & 3 & 0 & 0 \\
        $A_1$ & 22 & 24 & 3 & 4 & 2 \\
        $A_2$ & 25 & 24 & 3 & 7 & 4 \\
        $A_3$ & 29 & 24 & 3 & 12 & 5 \\
    \bottomrule 
    \end{tabular}
    \vspace{-2ex}
\end{table}

The benchmarks $A_1$, $A_2$, and $A_3$ mainly differ in the number of expansion options for the pipes and compressors. The expansion options for each of these networks is depicted in the Fig. \ref{fig:expansion-belgian}. For the location of the expansion options and the associated data including diameters, length, friction factors, etc. we use parameters from previous studies \cite{babonneau2012design, Borraz2016}. Though the compressor length and diameter are not explicitly optimized in any of the constraints involving the existing and new compressors, the data in these preceding studies specify those fields as well, and they are used in computing the the cost of construction of pipes and compressors, which is approximated by the following function \cite{babonneau2012design}: 
\begin{flalign}
c_p\!=\!\bm l_e ( 1.04081^{-6} \bm D_e^{2.5} \!+\! 11.2155) \,\, \forall (e, i, j) \!\in\! \mathcal C_n \cup \mathcal P_n. \label{eq:cost}
\end{flalign}
All the formulations for the RGEPP were implemented in the Julia programming language using optimization layer JuMP v0.18.6 \cite{Dunning2017} and GasModels v0.5.0\footnote{https://github.com/lanl-ansi/GasModels.jl}. Solver runs were performed on a MacBook Pro with a 2.9 GHz Intel Core i5 processor and 16GB RAM, with Gurobi v8.0 used to solve the MISOC relaxation of the RGEPP.

\subsection{Scenario Generation} \label{subsec:scenario-generation}
Here, we detail the method used to generate the uncertainty sets for each of the benchmarks. For each benchmark, a nominal set of gas demands at each receipt point in $\mathcal R$ is pre-specified as a part of the input data. The factors $\delta_1$, $\delta_2$, and $\delta_3$ are used to scale all the gas demands for the respective test cases $A_1$, $A_2$, and $A_3$. Based on these scaled gas demand values, we then create each scenario profile by defining a box uncertainty of width $2\varepsilon$, centered around the scaled gas withdrawal values for each receipt point. For the cases $A_1$ and $A_3$, we have a single profile, so $|\mathcal K| = 1$, whereas for $A_2$, we use two profiles and thus have $|\mathcal K| = 2$. Hence, for the case $A_2$, we have one scaling factor for each profile of scenarios $\mathcal S_2^1$ and $\mathcal S_2^2$, i.e., $\delta_2^{1}$ and $\delta_2^{2}$. The two scenarios for $A_2$ are meant to simulate a summer ($k=1$) and winter ($k=2$) cases. More scenarios for each of the benchmarks can be created, but we do not do so as it adds little qualitative value to the results. The scaling factors are set to the values $\delta_1 = 0.95$, $\delta_2^{1} = 1.0$, $\delta_2^{2} = 1.11$, and $\delta_3 = 1.0$ and the value of $\varepsilon$ is varied from $1\%$ to $5\%$ in steps of $1\%$. 

\subsection{Effectiveness of the MISOC Relaxation} \label{subsec:misoc-results}
The MISOC relaxation of the RGEPP is solved to optimality on all the benchmarks using the scenarios generated using the method in the previous section. The computation times for all the runs was less than a second and hence, we do not report any computation times for any of the runs in the case studies. The Tables \ref{tab:r-13} and \ref{tab:r-2} show the objective values of the expansion plans computed using the MISOC relaxation for the benchmarks $A_1$ and $A_3$, and $A_2$, respectively. It can be observed from both the tables that as we account for more uncertainty in the gas draws using the parameter $\varepsilon$, the cost of the resulting optimal expansion plan increases. 

\begin{table}[!h]
\vspace{-1ex}
    \centering
    \caption{Cost of expansion plans for the benchmarks $A_1$ and $A_3$}
    \label{tab:r-13}
    \begin{tabular}{ccc}
        \toprule 
        $\varepsilon$ & $A_1$ - objective value & $A_3$ objective value \\
        \midrule 
        \begin{filecontents*}{r1.csv}
epsilon,A1-obj,A2-obj,A3-s1-obj,A3-s2-obj,A3-obj
1,0.0,3206.59,1687.46,3409.59,3409.59
2,0.0,3206.59,1687.46,3409.59,3409.59
3,0.0,3206.59,1687.46,3409.59,3409.59
4,0.0,4987.2,3409.59,3409.59,3409.59
5,144.45,4987.2,3409.59,3409.59,3409.59
\end{filecontents*}
\csvreader[late after line=\\]{r1.csv}{1=\ep,2=\first,3=\third}{\ep & \first & \third}
        \bottomrule
    \end{tabular}
    \vspace{-1ex}
\end{table}

\begin{table}[!h]
    \vspace{-1ex}
    \centering
    \caption{Cost of expansion plans for the benchmark $A_2$}
    \label{tab:r-2}
    \begin{tabular}{cccc}
        \toprule 
        $\varepsilon$ & objective value  & objective value & objective value\\
        &   for $\mathcal S_2^1$ alone & for $\mathcal S_2^2$ alone & for both  scenarios \\
        \midrule 
        \csvreader[late after line=\\]{r1.csv}{1=\ep,4=\first,5=\second,6=\third}{\ep & \first & \second & \third}
        \bottomrule
    \end{tabular}
    \vspace{-1ex}
\end{table}

\subsection{Deterministic Network Expansion Planning Solution} \label{subsec:comparison}
\begin{figure}[!ht]
    \centering
    \includegraphics[width=\linewidth]{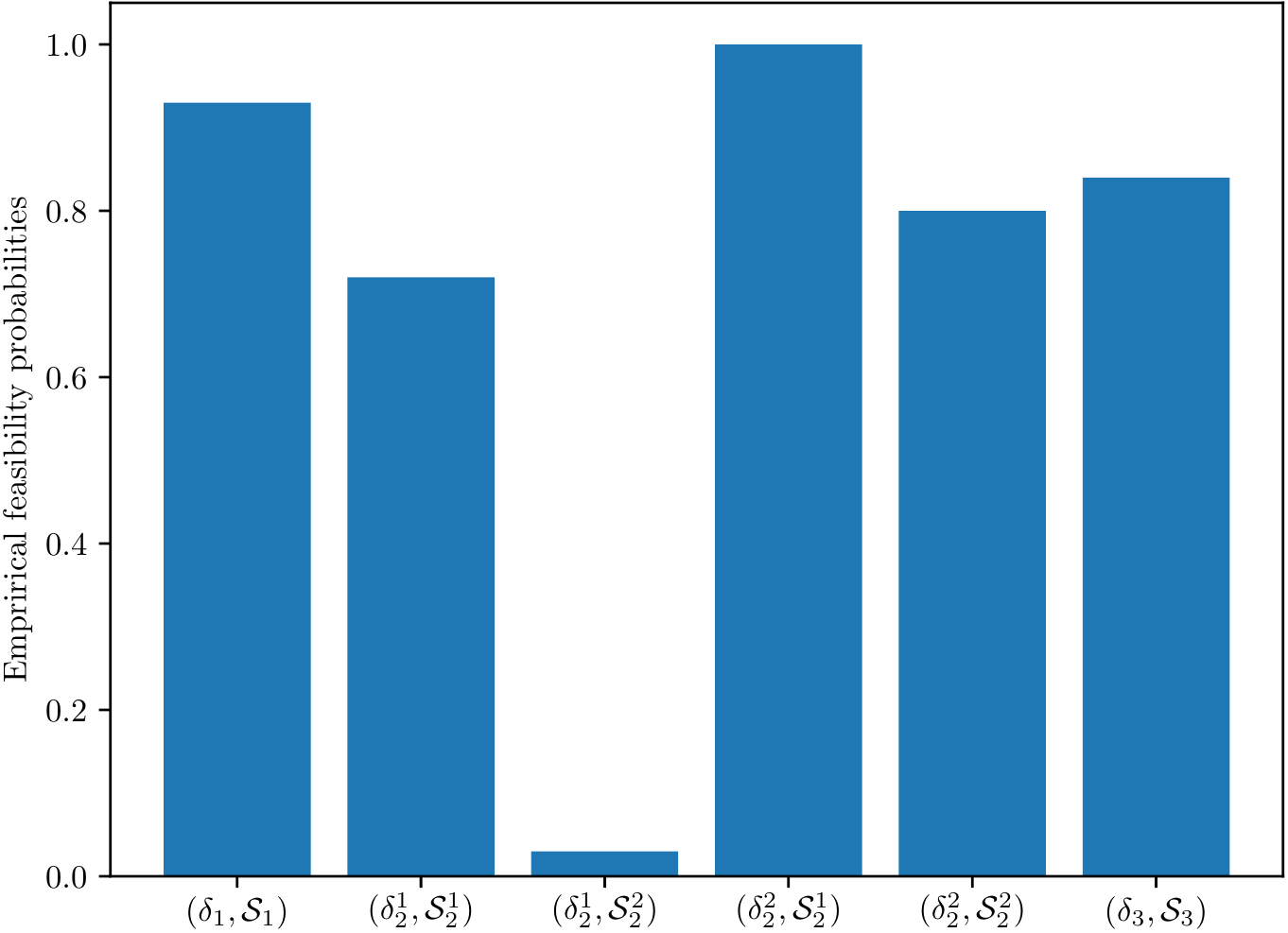}
    \caption{Empirical probabilities for randomly generated scenarios to be feasible for the expansion plan computed using nominal gas scenarios instead or uncertainty sets $\mathcal S$. The x-axis labels in this figure correspond (from left to right) to $(\delta_1,\mathcal S_1)$, $(\delta_2^1,\mathcal S_2^1)$, $(\delta_2^1,\mathcal S_2^2)$, $(\delta_2^2,\mathcal S_2^1)$, $(\delta_2^2,\mathcal S_2^2)$, and $(\delta_3,\mathcal S_3)$ in the notation used above.}
    \label{fig:prob}
\end{figure}

In this section, we compare a network expansion plan computed for a single (deterministic) scenario $\bm s$ that represents the mean demand with a robust expansion plan computed to account for the entire uncertainty set $\bigcup_{k\in\mathcal K}\mathcal S_k$. The deterministic network expansion planning problem yields a solution for one nominal profile unlike the RGEPP, which considers disjoint union of uncertain gas loads. The weakness of the deterministic formulation is shown by computing empirical probabilities of feasibility for the optimal plan computed using the deterministic formulation given a Monte Carlo sampling of the uncertain gas demand set. To compute these probabilities, we first solve the deterministic formulation for each benchmark using nominal gas withdrawal profiles computed using the scaling factors detailed in Sec. \ref{subsec:scenario-generation}. For the benchmarks $A_1$ and $A_3$, two deterministic runs (one for each benchmark) are performed using scaling factors $\delta_1$ and $\delta_3$, respectively. Then, $1000$ samples are generated from each of the scenario profiles $\mathcal S_1$ and $\mathcal S_3$ (computed using the parameter $\varepsilon$ set to $5\%$) and checked to see if they are feasible for the corresponding deterministic expansion plans. This feasibility result is used to compute empirical feasibility probabilities that are shown in Fig. \ref{fig:prob}. Similarly, for the benchmark $A_2$, we have two deterministic expansion plans for the scaling values $\delta_2^{1}$ and $\delta_2^{2}$; these expansion plans correspond to the nominal summer and winter gas loads. Each of these plans are checked for feasibility using samples from the summer and winter scenario profiles ($\mathcal S_2^{1}$ and $\mathcal S_2^{2}$, respectively). As seen from the figure, the network expansion plan obtained using the nominal winter demand (scaling $\delta_2^{2}$) for the benchmark $A_2^2$ is always feasible for any sample generated in $\mathcal S_2^{1}$, the summer scenario profile. The observation reverses for the pair $\delta_2^{1}$ and $\mathcal S_2^{2}$ where the probability is very low. Overall the RGEPP is always feasible for any sample generated from the uncertainty sets, whereas for the deterministic expansion plans there is always a finite probability that some scenario generated from an uncertainty set would be infeasible. 

\subsection{Effect of relaxing Monotonicity Constraints} \label{subsec:monotonicity-relaxation} 
In this section, we examine how relaxing the constraints that enforce monotonicity affect the network expansion plans computed using the MISOC relaxation of the RGEPP. We perform all the runs that were performed in Sec. \ref{subsec:misoc-results} and observed that despite relaxing the constraints \eqref{eq:compression_policy} that enforce the monotonicity-promoting compressor operation policy, the expansion plans produced by the resulting relaxed MISOC formulation remain exactly the same as those produced by the formulation that included these constraints. The only difference between the two solutions were the pressure profiles at the nodes i.e., the solutions with and without these constraints were considerably different in terms of the pressure profiles. 

\section{Conclusion and Future Work} \label{sec:conclusion}

In this study, we have examined the problem of optimal transport capacity expansion planning to enable a gas pipeline network to service the growing demand of gas-fired power plants.  Such generators are increasingly used to provide base load, flexibility, and reserve generation for bulk electric system. The aim of the developed robust gas pipeline network expansion planning problem (RGEPP) is to determine the minimal cost set of additional pipes and gas compressors that can be added to the network to provide the additional capacity to service future uncertain and variable loads.  We specifically seek to isolate the uncertainty that is caused by variability and intermittent nature of loads on a gas pipeline system that are caused by gas-fired generator commitment and economic dispatch, as determined by power system operations.  We do not examine uncertainties that may arise from resource allocation, changes in climate, and economics.

The RGEPP combinatorial optimization problem is initially formulated as a semi-infinite mixed-integer nonlinear program that accounts for box uncertainty constraints within multiple scenario profile sets. This formulation is then augmented by applying monotonicity properties of gas flows that promote desirable properties, which together with relaxation to a convex mixed integer second order cone problem (MISOCP) promotes computational tractability.  Our computational experiments support our claim for efficiency of the numerical method, and its ability to appropriately account for uncertainty in system loads. The resulting approach to solving a highly challenging semi-infinite mixed-integer nonlinear program enables expansion planning of gas pipeline systems to account for the variability that is inherent to the demands of gas-fired electricity production.  While we consider steady-state flow modeling here, a promising future direction involves extending the formulation so that the modeled capacity estimate and resulting expansion plan would truly account for the intra-day intertemporal variability of gas-fired generator fuel usage throughout an integrated natural gas and electric energy delivery system.

\section*{Acknowledgements}
This work was supported by the U.S. Department of Energy's Advanced Grid Modeling (AGM) project \emph{Joint Power System and Natural Gas Pipeline Optimal Expansion}. The research work conducted at Los Alamos National Laboratory is done under the auspices of the National Nuclear Security Administration of the U.S. Department of Energy under Contract No. 89233218CNA000001. We gratefully thank the AGM program manager Alireza Ghassemian for his support of this research.

\bibliographystyle{IEEEtran}
\bibliography{expansion}

\end{document}